\documentclass[10pt]{article}

\usepackage[utf8]{inputenc}
\usepackage[T1]{fontenc}
\usepackage[english]{babel}
\usepackage{cmap}

\usepackage{graphicx}
\usepackage{amsmath}
\usepackage{amssymb}
\usepackage{amsthm}

\newcommand{\X}{\mathfrak X}

\newcommand{\romb}{\diamondsuit}
\newcommand{\Log}{{Log}}

\newcommand{\logic}[1]{\mathsf{#1}}

\newcommand{\lK}{\logic{K}}

\newcommand{\Lo}{\mathsf{L}}

\newcommand{\pmor}{\twoheadrightarrow}

\newcommand{\set}[1]{\left\{#1\right\}}

\newcommand{\setdef}[2][x]{\set{#1\,\left |\,#2 \right .}}

\newcommand{\pair}[1]{\left \langle #1 \right \rangle}

\newcommand{\QQ}{\mathbb Q}
\newcommand{\Natr}{\mathbb N}

\newcommand{\eps}{\varepsilon}

\newcommand{\Y}{\mathcal Y}

\newcommand{\st}{\mathop{st}}
\newcommand{\Top}{\mathop{Top}}

\renewcommand{\root}{\circledR}

\newtheorem{theorem}{\textsc{Theorem}}[section]
\newtheorem{corollary}[theorem]{\textsc{Corollary}}
\newtheorem{lemma}[theorem]{\textsc{Lemma}}
\newtheorem{proposition}[theorem]{\textsc{Proposition}}

\theoremstyle{remark}

\newtheorem*{remark*}{\textmd{\textsc{Remark}}}

\theoremstyle{definition}
\newtheorem{definition}{{Definition}}[section]
\newtheorem*{definition*}{{Definition}}

\begin{document}
	
		\title{Topological product of modal logics with the McKinsey axiom}
		\author{Andrey Kudinov}
		
		\maketitle
		
		Steklov Mathematical Institute of RAS\footnote{
		The research was supported by RSF (project No. 21-11-00318)}
		
		\begin{abstract}
			In this paper we consider the topological products of modal logics of S4.1 and S4.  We prove that it is equal to the fusion of logics S4.1 and S4 plus one additional axiom. We also show that this product is decidable. This is an example of a topological product of logics that is greater than the fusion but less than the expanding product of the corresponding logics.
		\end{abstract}
		

	\section{Introduction}
	Topological products of modal logics were introduced in 2006 in the work of J. van Bentham and co-authors \cite{benthem:MultimodalLogicsProductsTopologies}. 
	The product of two modal logics is the logic of the class of all products of structures of the corresponding logics.
	The topological product of modal logics $\logic{L_1}$ and $\logic{L_2}$ containing $\logic{S4}$ is defined semantically as the modal logic of the class of all possible products of topological spaces of corresponding logics. 
	Note that the product of topological spaces defined in \cite{benthem:MultimodalLogicsProductsTopologies} differs from the product known to us from topology. 
	To distinguish between these two concepts, we will call the van Bentham construction the \emph{bitopological product}, because the result is a space with two topologies. In \cite{benthem:MultimodalLogicsProductsTopologies}, the topological product $\logic{S4}\times_t\logic{S4}$ was shown to coincide with the fusion $\logic{S4}\ast\logic{S4}$. Completeness with respect to the bitopological product $\QQ \times\QQ$ was also established in \cite{benthem:MultimodalLogicsProductsTopologies}.

	In \cite{kremertopological}, Ph.~Kremer proved that the topological product of the logics $\logic{S4}$ and $\logic{S5}$ is equal to the semi-commutator (or expanding commutator) of these logics.
	$$ 
	[\logic{S4}, \logic{S5}]^{EX} = \logic{S4} \ast \logic{S5} + \Box_1 \Box_2 p \to \Box_2 \Box_1 p + \romb_1\Box_2 p \to \Box_2\romb_1 p.
	$$ 
	According to the well-known completeness result for expanding products of Kripke frames \cite[\S 9.1]{gabbay:Manydimensional} this logic coincides with the expanding product of these logics $(\logic{S4} \times \logic{S5})^{EX}$.
	
	In the paper by Kremer \cite{kremertopological}, it is shown that the topological product of two modal logics $\Lo_1$ and $\Lo_2$ must be between their fusion $\Lo_1\ast\Lo_2$ and the product $\Lo_1\times\Lo_2$. So far, we know a topological product of logics that is equal to the fusion, namely $\logic{S4} \times_t\logic{S4}$ (see \cite{benthem:MultimodalLogicsProductsTopologies}). There is a topological product that is equal to the usual product ($\logic{S5} \times\logic{S5} $, see \cite{kremertopological}) and the product that is equal to the expanding product ($\logic{S4} \times_t \logic{S5}$, see \cite{kremertopological}). Other topological products of modal logics expending $\logic{S4}$ have not been found yet.
	
	In this paper, we consider the topological product of the logics $\logic{S4.1}$ and $\logic{S4}$ and prove that it is equal to $\logic{S4.1}\ast\logic{S4} + \romb_1\Box_2(\romb_1 p\to \Box_1 p) $ and is strictly between the fusion of the logics $\logic{S4.1} \ast \logic{S4}$ and the expanding product of these logics. This is the first example of such a product
	
	\section{Definitions and known results}
	
	Let $ \mathrm{PROP} $ be a countably infinite set of propositional letters then we use the Backus-Naur form to define a \emph{modal formula} inductively:
	$$
	A ::= p\; |\;\bot \; | \; (A \to A) \; | \; \Box_i A,
	$$
	where $p\in \mathrm{PROP}$ is a propositional letter, and $\Box_i$ is a modal operator ($ i= 1, \ldots, N$).
	Other connectives are introduced as abbreviations: classical connectives are expressed through $\bot$ and $\to$,
	and $\Diamond_i$ is defined as $\lnot\Box_i\lnot$.
	In this paper we use modal languages with $ N \in \set{1,2} $. If there is only one modality, we omit the index.
	
	\begin{definition} A \emph{(normal) modal logic\/} is
		a set of modal formulas closed under Substitution $\left
		(\frac{A(p)}{A(B)}\right )$, Modus Ponens $\left (\frac{A,\,
			A\to B}{B}\right )$ and Necessitation $\left
		(\frac{A}{\Box_i A}\right)$ rules, containing all the
		classical tautologies and the normality axioms:
		$$
	\begin{array}{l}
		\Box_i (p\to q)\to (\Box_i p\to \Box_i q).
	\end{array}
	$$
	
	$\logic{K_N}$ denotes \emph{the smallest normal modal logic with $N$ modalities} and we write $\lK$ instead of  $\logic{K_1}$.
\end{definition}

Let $\logic{L}$ be a logic and $\Gamma$ be a set of formulas. Then $\logic{L} +
\Gamma$ denotes the smallest logic containing $\logic{L}$ and $\Gamma$. For
$\Gamma = \set{A}$ we write $\logic{L} + A$ rather than $\logic{L}+  \set{A}$.

\begin{definition}
	Let $\logic{L_1}$ and $\logic{L_2}$ be two modal logics with one modality  $\Box$ (unimodal logics), then the \emph{fusion} of these logics is the following modal logic with 2 modalities:
	\[
	\logic{L_1} \ast \logic{L_2} = \logic{K_2} + \logic{L'_1} + \logic{L'_2};
	\] 
	where $\logic{L'_i}$ is the set of all formulas from $\logic{L_i}$ with all instances of $\Box$ being replaced with $\Box_i$.
\end{definition}

Let us define logic $\logic{S4}$ in a standard way:
\[ 
\logic{S4} = \logic{K} + \Box p \to p + \Box p \to \Box \Box p.
\]

\begin{definition}
	A \emph{Kripke $N$-frame} is a tuple $ F = (W, R_1, \ldots, R_N)$, where $W\ne \varnothing$ is a set. Elements of $W$ we call \emph{possible worlds} or \emph{point} and $R_1, \ldots, R_N \subseteq W \times W$ are relations on $W$.
\end{definition}

\begin{definition}
	A \emph{valuation} on a Kripke frame $F = (W, R_1, \ldots, R_N)$ is a function $V:\mathrm{PROP} \to \mathcal{P}(W)$, where $\mathcal{P}(W)$ is the set of all subsets of $W$. Pair $M = (F, V)$ is called a \emph{Kripke model}. The truth relation ``$\models$'' at a point in a model is defined recursively:
	\begin{align*}
		M, x \models p &\iff x \in V(p),\ \hbox{for $p\in \mathrm{PROP}$};\\
		M, x \not\models \bot; &\\
		M, x \models A \to B &\iff \bigl(M, x \models A \Rightarrow M, x \models B \bigr);\\
		M, x \models \Box_i A &\iff \forall y (xR_i y \Rightarrow M,y \models A).
	\end{align*}

	Formula $ A $ is \emph{valid} in a frame $ F $ if $ \forall V \forall x\in W (F,V,x \models A)$ (Notaion: $ F \models A $).
\end{definition}

\begin{definition}
	Let $ X \ne \varnothing $. A \emph{topology}  on $ X $ is a collection $ T $ of subsets of $ X $ satisfying the following conditions:
	\begin{itemize}
		\item $\varnothing, X \in  T$;
		\item if $ U_1, U_2 \in T $ then $ U_1 \cap U_2 \in T $;
		\item if $ S \subseteq T $ then $ \bigcup S \in T $.
	\end{itemize}
	Pair $ (X, T) $ is called a \emph{topological space}. Elements of $T$ are called \emph{open sets}.
	
	A collection of open subsets $B$ of a topological space is called \emph{a base for the topology} if every open subset is a union of some elements of $B$.
	
	\emph{A bitopological space} is a triple $(X, T_1, T_2)$, where $T_1$ and $T_2$ are topologies on $X$.
\end{definition}

\begin{definition}
	A \emph{valuation} on a topological space $\X = (X, T)$ is a function $V:\mathrm{PROP} \to \mathcal{P}(X)$. Pair $M = (\X, V)$ is called a \emph{topological model}. The truth relation is defined as in Kripke model with a difference in the last point:  
	\begin{align*}
		\mathcal{M}, x \models \Box A &\iff \exists U \in T \left ( x \in U \; \& \; \forall y \in U (\mathcal{M},y \models A) \right ).
	\end{align*}
	
	This definition can be naturally extended to the language with two modalities. In this case, the models are based on bitopological spaces.
	
	Formula $ A $ is \emph{valid} in a space $ \X $ if $ \forall V \forall x\in W (\X,V,x \models A)$ (notation: $ \X \models A $).
\end{definition}

For a class of topological spaces (Kripke frames) $ \mathcal{C} $ we define \emph{the logic of $ \mathcal{C} $} as
\[  
Log(\mathcal{C}) = \setdef[A]{\forall S \in \mathcal{C} (S\models A)}.
\] 
For a one-element class we omit the curly braces and write $Log(\X)$. 

Let $ F = (W, R) $ be a Kripke frame and $ F \models \logic{S4} $ (or, equivalently, $ R $ is transitive and reflexive, see \cite[\S 3.8]{Chagrov&Z97}). We define topological space $\Top(F) = (W, T_R)  $, where $ T_R $ is a topology with base $ \setdef[R(x)]{x\in W} $.

\begin{definition}
	Let $F = (W, R_1, \ldots, R_N)$ and $G = (U,S_1, \ldots, S_N)$ be two Kripke frames. Function $f:W \to U$ is a \emph{p-morphism}  from $F$ onto $G$ if
	\begin{enumerate}
		\item $f$ is surjective;
		\item $xR_iy \Rightarrow f(x) S_i f(y)$ (monotonicity);
		\item $f(x) S_i u \Rightarrow \exists y (f(y) = u \;\&\; xR_iy)$ (lifting).
	\end{enumerate}
	
	Notation: $f: F \twoheadrightarrow G$.
\end{definition}

The following theorem is well-known (see \cite{Chagrov&Z97}).
\begin{theorem}[p-morphism]\label{thm:pmorphism}
	If $F \twoheadrightarrow G$ then $Log(F) \subseteq Log(G).$
\end{theorem}

For a reflexive and transitive frame $F$ we can define $Log(F) = Log(Top(F))$. In general $Log(F)$ is the set of all valid formulas in $F$.

\begin{definition}
	Let $\X$ and $\Y$ be two topological spaces. Function $f:\X \to \Y$ is  called \emph{open} (\emph{continuous}) if the image (preimage) of any open set is open.
\end{definition}

The analogue of a p-morphism for topological spaces is a surjective, open and continuous map. The corresponding theorem is also true (see \cite{gabelaia_modal_2001}): 
\begin{theorem}[topological p-mophism]\label{thm:top_pmorphism}
	Let $\X$ and $\Y$ be topological spaces and $f: \X\to \Y$ be a surjective, open and continuous map. Then $Log(\X) \subseteq Log(\Y).$
\end{theorem}

Such maps we will also call \emph{p-morphisms}. There will be no collisions because of the following
\begin{lemma}
	For Kripke frames $ F_i = (W_i, R_i) $, such that $ R_i $ is reflexive and trasitive ($ i = 1,2 $) a surjective map $ f:W_1 \to W_2 $ is a p-morphism iff $ f $ is open and continuous with respect to $ T_{R_1}  $ and $ T_{R_2}$ topologies.
\end{lemma}
The proof is rather straightforward and we leave it to the reader.

\begin{definition}
	Let $ F = (W, R_1, R_2) $ be a 2-frame such that both $R_1$ and $R_2$ are reflexive and transitive. We define a bitopological space $\Top_2(F) = (W, T_{R_1}, T_{R_2})  $.
\end{definition}

\begin{definition}(\cite{benthem:MultimodalLogicsProductsTopologies})
	Let $ \X_1 = (X_1, T_1)$ and $ \X_2 = (X_2, T_2) $ be two topological spaces. We define the (bitopological) product of them as the bitopological space $ \X_1 \times \X_2 = (X_1 \times X_2, T_1^h, T_2^v )$. Topology $ T_1^h $ is the topology with the base $ \setdef[U\times \set{x_2}]{U \in T_1\;\&\; x_2\in X_2} $ and topology $ T_2^v $ is the topology with the base $ \setdef[\set{x_1}\times U]{x_1\in X_1 \;\&\;  U \in T_2} $.  We call topolgy $ T_1^h $ \emph{horisontal} and topology $ T_2^v $ \emph{vertical}.
\end{definition}

\emph{The topological product of 
modal logics} $ L_1 $ and $ L_2 $ is the following logic with two modalities:
\[  
L_1 \times_t L_2 = Log(\setdef[\X_1 \times \X_2]{ \X_1, \X_2\hbox{ --- topological spaces,}\ \X_1 \models L_1, \X_2 \models L_2}).
\]

\begin{theorem}[\cite{vanBenthem2007modal}]
	$\logic{S4} \times_t \logic{S4} = \logic{S4} \ast \logic{S4}.$
\end{theorem}

\section{McKinsey axiom}

Formula $ A1 = \Box \romb p \to \romb \Box p $ is known in the literature as the \emph{McKinsey axiom}. It is well-studied in the context of Kripke semantics. 
In \cite{Goldblatt91} it was shown that this formula is not canonical\footnote{The definition of a canonical formula is given at the beginning of Section \ref{sec:MainResults}.}. The topological aspects of this axiom was studied in \cite{vanBenthem2007modal,bezhanishvili2012modal,bezhanishvili2015modal}).

\begin{lemma}\label{lem:McKinseyProperty}
For a transitive frame $ F = (W,R) $ the validity of formula $A1$ is equivalent to the following first-order property 
\[ 
\forall w \in W \exists u \in W (w R u \land R(u) = \set{u} ),
\]
where $ R(u) = \setdef[t]{uRt} $.
\end{lemma} 
The proof can be found in \cite{Chagrov&Z97}.

\begin{definition}
	In a topological space $\X$ point $ x $ is \emph{isolated} if set $ \set{x} $ is open. 
	$ \X $ is
	\emph{weakly scattered} if the set of all isolated points of $ \X $ is dense in $ \X $, that is if any open subset includes an isolated point.
\end{definition}

\begin{theorem}[\cite{gabelaia_modal_2001}]
	$ \logic{S4.1} $ is the logic of the class of weakly scattered spaces.
\end{theorem}

\begin{lemma}
	Let $ \X_1 $ and $ \X_2 $ be topological spaces and $ \X_1 $ is weakly scattered. Then
	\[ 
	\X_1 \times \X_2 \models \romb_1\Box_2(\romb_1 p \to \Box_1 p).
	\]
\end{lemma}

\begin{proof}
	Let us take $ \pair{x,y} \in \X_1 \times \X_2 $ and a horizontal open neighborhood $ U\times \set{y} $, where $ U \in T_1 $ and $ \X_1 = (X_1,T_1) $. Since $ \X_1 $ is weakly scattered, set $ U $ contains an isolated (in $ \X_1 $) point $ x' $.
	
	It follows that for any $ y'\in \X_2 $ point $ \pair{x',y'} $ is isolated in horizontal topology, and hence
	$\pair{x',y'} \models \romb_1 p \to \Box_1 p$. Then 
	$ 
	\pair{x',y} \models \Box_2(\romb_1 p \to \Box_1 p)
	$
and  $\pair{x,y} \models \romb_1\Box_2(\romb_1 p \to \Box_1 p).$
	
\end{proof}

\begin{lemma}\label{lem:MKDuo_prop}
	For any $ \logic{S4.1} \ast \logic{S4}$-frame $ F = (W, R_1, R_2) $ it is true that
	\[ 
	F \models \romb_1\Box_2(\romb_1 p \to \Box_1 p) \iff \forall x \exists y (x R_1 y \ \& \ \forall z (y R_2 z \Rightarrow R_1 (z) = \set{z})).
	\]
\end{lemma}
The proof can be easily obtained by modifying the proof of Lemma   \ref{lem:McKinseyProperty}. We leave the details to the reader.

\begin{lemma}
	$ \romb_1\Box_2(\romb_1 p \to \Box_1 p) \notin \logic{S4.1} \ast \logic{S4}$.
\end{lemma}

\begin{figure}[h]
	\centering
	\includegraphics[width=0.4\linewidth]{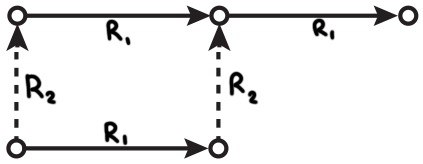}
	\caption{}
	\label{fig:frame1}
\end{figure}

This lemma can be proved by providing a $ \logic{S4.1} \ast \logic{S4} $-frame that falsify the formula. Please see such a frame on Fig.\ref{fig:frame1}.

\begin{corollary}
	$\logic{S4.1} \times_t \logic{S4} \ne \logic{S4.1} \ast \logic{S4}$.
\end{corollary}

\begin{lemma}\label{lem:MKDuo_prop}
	$ \romb_1\Box_2(\romb_1 p \to \Box_1 p) \in (\logic{S4.1} \times \logic{S4})^{EX}$.
\end{lemma}

The proof is left to the reader.

\section{Completeness and decidability theorems}
The canonical model construction is a well-known method that is often used to prove the Kripke completeness of a modal logic. A modal logic is called \emph{canonical} if all of its formulas are valid in its corresponding canonical frame. The canonical frame of any logic falsify all formulas that are not in that logic. Therefore, any canonical logic is also Kripke-complete.
For details see \cite{Chagrov&Z97,blackburn_modal_2002,NonclassicalLogics_Odintsov_2014_eng}.

Let us denote
$$ 
\Lo = \logic{S4.1} \ast \logic{S4} + \romb_1\Box_2(\romb_1 p \to \Box_1 p).
$$

\begin{theorem}
	Logic $ \Lo $ is canonical and, as a corollary,  Kripke complete.
\end{theorem}

Indeed, logic $ \logic{S4.1} \ast \logic{S4} $ is canonical since all its axioms are canonical. The canonicity of $ \romb_1\Box_2(\romb_1 p \to \Box_1 p) $ can be obtatianed by a straightforward modification of the proof of canonisity of McKinsey axiom (see \cite[Ch. 5, Th. 5.21]{Chagrov&Z97}).
The details of the proof are easy to reconstruct.

A logic is said to have the \emph{finite model property} if it is the logic of a class of finite frames\footnote{The finite model property is equivalent to the finite frame property (cf. \cite[Section 8.4]{Chagrov&Z97}).}.

\begin{theorem}[Harrop, cf. {\cite[\S 16.2]{Chagrov&Z97}}]
	Any finitely axiomatizable logic which has finite model property is decidable. 
\end{theorem}

\begin{theorem}
	Logic $ \Lo $ has the finite model property and is decidable.
\end{theorem}

\begin{proof}
To prove the finite model property of $ \Lo $ we use filtration with prior partitioning. The details of this method can be found in \cite{KuSh2012}. Let $ A $ be a formula and $ M = (F,V) $ a model such that $ F\models \Lo $, $ M \not\models A $ and $ F = (W, R_1, R_2) $. We define partition $ W $ into to two subsets
$  
W = W_1 \sqcup W_2,
$
where $ W_1 $ is the set of all $ R_1 $-maximal points and $ W_2 = W \setminus W_1 $. 
Let us define relation $ \sim $:
\begin{align*}
x \equiv_A y &\iff \forall B (B \ \hbox{is a subformula of $ A $ and}\ (M,x \models B \Leftrightarrow M,y \models B));\\
x \sim y &\iff x \equiv_A y \ \hbox{and}\ \exists i \in \set{1,2} (x,y \in W_i).
\end{align*}

Let $ M' = (F', V')$ be the transitive filtration of $ M $ via $ \sim $. It is easy to check that $ R_1 $-maximal points will be preserved and hence $ F' \models \Lo $. By the filtration lemma we have $ M' \not\models A $ and $ M' $ is finite. The finite model property is proven. 
\end{proof}

The main result of this paper is the following theorem.

\begin{theorem}\label{thm:main}
	$\logic{S4.1} \times_t \logic{S4} = \logic{L}.$
\end{theorem}

Let $ \mathbb{T}_{2,2} $ be the infinite transitive $ (2,2) $-tree with two relations:  $ \mathbb{T}_{2,2} = (T_{2,2}, R_1, R_2) $, where $T_{2,2} = \set{a_1, a_2, b_1, b_2}^*$ is the set of finite words in a 4-letter alphabet. For any $ \vec a, \vec b \in  T_{2,2}$ we put
\begin{align*}
	\vec a R_1 \vec b &\Leftrightarrow \hbox{there exists $ \vec c\in \set{a_1, a_2}^* $, such that} \ \vec b  = \vec a \cdot \vec c;\\
	\vec a R_2 \vec b &\Leftrightarrow \hbox{there exists $ \vec d\in \set{b_1, b_2}^*$, such that}\ \vec b = \vec b \cdot \vec d.
\end{align*} 
Here and further on $\cdot$ denote the concatenation operation on words. 

By $ \mathbb{T}_{2,2+2} = (W, R_1', R_2')$ we denote a frame that we get by putting a copy of the infinite reflexive and transitive 2-tree $ \mathbb{T}_2 = (T_2, \sqsubseteq) $ above every point in $ \mathbb{T}_{2,2} $, here $T_2 = \set{1,2}^*$ and $\sqsubseteq$ is the relation of being a prefix.  
Let $ \eps $ be the empty word which is the root of $ \mathbb{T}_2 $. The copies of $\mathbb{T}_2$ themselves are connected by $ R'_1 $, while inside the trees the points are connected by $R'_2$. The precise definition is as follows:
\begin{enumerate}
	\item\label{it:i1}  $ W = T_{2,2} \times \set{\root} \cup T_{2,2} \times T_2$, where $\root$ is a new special symbol;
	\item\label{it:i2} 
	$ \vec a R_1 \vec a' \Rightarrow \pair{\vec a, \root} R'_1 \pair{\vec a', \root}$;
	\item\label{it:i3} 
	$ \pair{\vec a, \root} R'_1 \pair{\vec a, \eps} $ for every $\vec a$;
	\item $R'_1$ is the minimal reflexive and transitive relation which satisfies items (\ref{it:i2}) and (\ref{it:i3});
	\item\label{it:i4}
	$ \vec a R_2 \vec a' \Rightarrow \pair{\vec a, \root} R'_2 \pair{\vec a', \root}$;
	\item\label{it:i5}  $\vec{a} = \vec{a}'$, $ \vec b \sqsubseteq \vec b' $ and $ \vec b,\vec b' \in T_2 \Rightarrow \pair{\vec a, \vec b} R'_2 \pair{\vec a', \vec b'}$ ,
	\item $R'_2$ is the minimal reflexive and transitive relation which satisfies items (\ref{it:i4}) and (\ref{it:i5}).
\end{enumerate}

From this definition it follows 
\begin{lemma}The following equities are true
	\begin{align*}
		R'_1(\pair{\vec a, \root}) &= R_1(\vec a) \times \set{\root, \eps};\\ 
		R'_1(\pair{\vec a, \vec b}) &= \set{\pair{\vec a, \vec b}},\hbox{ for } \vec b \in T_2;\\ 
		R'_2(\pair{\vec a, \root}) &= R_2(\vec a) \times \set{\root};\\ 
		R'_2(\pair{\vec a, \vec b}) &= \set{\vec a} \times {\sqsubseteq}(\vec b),\hbox{ for } \vec b \in T_2.
	\end{align*}
\end{lemma}

\begin{lemma}
	$ \Log(\mathbb{T}_{2,2+2}) = \Lo$.
\end{lemma}

\begin{proof}
	It is easy to check that $ \mathbb{T}_{2,2+2} \models \Lo $.
	
	Let us assume that $A \notin \Lo$. Logic $\Lo$ has the finite model property, so there exists a rooted finite frame $F$, such that $F \models \Lo$ and $F \not\models A$. To finish the proof it is sufficient to show that
	$ \mathbb{T}_{2,2+2} \pmor F $.
	
	In \cite{vanBenthem2007modal} a p-morphism
	$ f: \mathbb{T}_{2,2} \pmor F$ was described.

	By Lemma \ref{lem:MKDuo_prop} for every $ w\in F $ there exists $ \mu(w) \in R_1(w) $ such that all $ u \in R_2(\mu(w)) $ are $ R_1$-maximal. And for each $ w$ we fix a p-morphism $ h_w: \mathbb{T}_2 \pmor \left (R_2(\mu(w)), R_2|_{R_2(\mu(w))}\right ) $. Now we define $ g: \mathbb{T}_{2,2+2} \pmor F $:
	\begin{align*}
	g\left (\pair{\vec a, \root}\right ) &= f(\vec a);\\
	g\left (\pair{\vec a, \vec b}\right ) &=h_{f(\vec a)}(\vec b).
	\end{align*}

Let us check that $ g $ is indeed a p-morphism.

The surjectivity of $g$ follows from the surjectivity of $f$. 
The monotonicity follows from the monotonicity of $f$ and $h_w$ and transitivity of all relations. 
The lifting follows from the lifting of $f$ and $h_w$. 
It is sufficient to consider two cases: $x = \pair{\vec a, \root}$ and $ x = \pair{\vec a, \vec b}$.
\end{proof}

To finish the proof of Theorem \ref{thm:main} we need to construct a weakly scattered space $\X$ and a space $ \Y $ such that there exists an open and continuous surjection $ f: \X \times \Y \to Top_2(\mathbb{T}_{2,2+2})$.
In order to define such $ \X $ and $ \Y $ we first introduce several new notions.

\begin{definition}
	\emph{A path with stops} on $\mathbb{T}_2 = (\set{1,2}^*, \sqsubseteq)$ is a tuple $x_1 \ldots x_n$, where $x_i \in \set{0,1,2}$. We define $ f_F $ recursively on the set of all paths with stops:
	\begin{itemize}
		\item $ f_F(\eps) = \eps $;
		\item $ f_F(\vec a 0) = f_F(\vec a) $;
		\item $ f_F(\vec a 1) = f_F(\vec a)1 $;
		\item $ f_F(\vec a 2) = f_F(\vec a)2 $.
	\end{itemize}

\end{definition}

\begin{definition}
	A \emph{pseudo-infinite path} on $\mathbb{T}_2$ is an infinite sequence of 0, 1 and 2 which contains only finitely many non-zeros numbers. The infinite sequence of zeros is denoted as $0^\omega$. So $\alpha$ is a pseudo-infinite path if $\alpha = \vec{a}0^\omega$ for some $\vec{a} \in \set{0,1,2}^*$.
	Let $W_\omega$ is the set of all pseudo-infinite paths on $\mathbb{T}_2$. Next, we will define functions $st: W_\omega \to \Natr$ and $f_\omega:W_\omega \to T_2$.
	Let $\alpha = x_1 \ldots x_n \ldots \in W_\omega$, then
	\begin{align*}
		\st(\alpha) &= \min \setdef[N]{\forall k> N (a_k = 0)};\\
		\alpha\lceil_k &= x_1 \ldots x_k;\\ 
		f_\omega(\alpha) &= f_F\left (\alpha\lceil_{\st(\alpha)}\right );\\
		U_k(\alpha) &= \setdef[\beta \in W_\omega]{\alpha \lceil_k = \beta \lceil_k \ \&\ f_F(\alpha\lceil_k) R f_\omega(\beta)},\hbox{ where }\alpha \in W_\omega, k\in \Natr. 		
	\end{align*}
\end{definition}

\begin{lemma}
	A family of sets $B = \setdef[U_{k}(\alpha)]{\alpha \in W_\omega, k>0}$ is a topological base.
\end{lemma}
\begin{proof}
	Note that for arbitrary $\alpha, \beta \in W_\omega$ and $0<k\le m$ one of the following two items is true: 
	\begin{enumerate}
		\item $U_k(\alpha) \cap U_m(\beta) = \varnothing$;
		\item $U_k(\alpha) \subseteq U_m(\beta)$.
	\end{enumerate}
	Hence, $B$ is closed under intersections and is a topological base.
\end{proof}

Let $T_\omega$ be the topology generated by base $B$ and $\Y = (W_\omega, T_\omega)$. 

Next, we will define function $g:W_\omega\times W_\omega \to T_{2,2}$. Let $\alpha = x_1x_2\ldots \in W_\omega$ and $\beta = y_1y_2\ldots \in W_\omega$. Let us remind that $T_{2,2} = \set{a_1, a_2, b_1, b_2}^*$. For convenience, we assume that $a_0 = b_0 = 0$.
\[  
g(\alpha, \beta) = f'_\omega(a_{x_1}b_{x_1}a_{x_2}b_{x_2}\ldots),
\]
where $f'_\omega$ is the function that analogous to $f_\omega$ and deletes zeros in infinite sequences of symbols from set $\set{0, a_1, a_2, b_1, b_2}$. Since $\alpha$ and $\beta$ have tails of zeros, then $g(\alpha, \beta)$ be a finite sequences, i.e. an element of $T_{2,2}$.

\begin{proposition}[\cite{kudinov_aiml12}]\label{prop:pmorphismg}
	Function $g:\Y\times \Y \pmor Top_2(\mathbb{T}_{2,2})$ is a p-morphism.
\end{proposition}

Let us denote set $\setdef[k\in \Natr]{k\ge n}$ as $\Natr_{\ge n}$.

We define $\X = (X, T)$ as follows:
\begin{align*}
	X &= W_\omega \times \Natr,\\	
	U'_k (\alpha,0) &= \left (U_k(\alpha) \times \set{0}\right ) \cup \left (U_k(\alpha) \times \Natr_{\ge k}\right ), \\
	U'_k (\alpha,n) &= \set{\pair{\alpha,n}}, \hbox{ where } n\ge 1. 
\end{align*}
Family of sets $U'_k(\alpha, n)$ forms a base for topology $T$. 
To check the correctness of this definition it is sufficient to show that any two sets of type $U'_k(\alpha, n)$ are either do not intersect, or contain each other.

Points $\pair{\alpha, n}$  ($n\ge 1$) are isolated, hence any element of the base contains an isolated point. It follows 
\begin{lemma}\label{lem:S41xS4corr}
	Topological space $\X$  is weakly scattered and $\X \models A1$.
\end{lemma}

\begin{lemma}\label{lem:S41xS4comp}
	$\X \times \Y \pmor \Top_2 (\mathbb{T}_{2, 2+2})$.
\end{lemma}

\begin{proof}
	To construct a needed p-morphism we will use the p-morphism $g:\Y\times \Y \pmor Top_2(\mathbb{T}_{2,2})$ from Proposition \ref{prop:pmorphismg}. 
	
We define $ f:\X \times \Y \to \mathbb{T}_{2,2 + 2} $ as follows:	\begin{align*}
		f(\pair{\alpha, 0}, \beta) &= \pair{g(\alpha, \beta), \root} ;\\
		f(\pair{\alpha, n}, \beta) 
		&=\pair {g(\alpha\lceil_n \cdot 0^\omega, \beta\lceil_n \cdot 0^\omega), f_{\omega}(\gamma)},\hbox{ where } \beta = \beta\lceil_n\;\cdot\;\gamma.
	\end{align*}
	Surjectivity of $f$ follows from surjectivity of $g$ and $f_\omega$.

Let us check that $f$ is open and continuous.

($T^h_1$-openness) Let us take a horizontal neighborhood $
U'_k(\alpha, 0) \times \set{\beta} 
$ of point $(\pair{\alpha, 0}, \beta)$.
We can assume that $ k \ge \max(\st(\alpha), \st(\beta)) $ since any set $U'_k(\alpha, 0)$ is a union of similar sets but with greater $k$. The image of this neighborhood equals 
\begin{align*}
	&f(U'_k(\alpha, 0) \times \set{\beta}) = \\
	&f\left ( (U_k(\alpha) \times \set{0})\times \set{\beta} \cup \left (U_k(\alpha) \times \Natr_{\ge k}\right  ) \times \set{\beta}\right )=\\
	&f\left ( (U_k(\alpha) \times \set{0})\times \set{\beta}\right ) \cup f\left (\left (U_k(\alpha) \times \Natr_{\ge k}\right  ) \times \set{\beta}\right )=\\
	&g(U_k(\alpha) \times \set{\beta}) \times \set{\root, \eps}.
\end{align*}
The last equality is true since if $\beta = x_1x_2\ldots x_n x_{n+1}\ldots $ and $ n\ge k $. So
\begin{align*}
	f(\pair{\alpha, n}, \beta) = &\pair{g(\alpha\lceil_n \;\cdot\; 0^\omega, \beta\lceil_n \;\cdot\; 0^\omega), f_{\omega}(x_{n+1} x_{n+2}\ldots)} =\\
	&\pair{g(\alpha, \beta), f_{\omega}(0^\omega)} = \pair{g(\alpha, \beta), \varepsilon}.
\end{align*}
 
 Point $ \pair{g(\alpha, \beta), \varepsilon} $ is $T_1^h$-isolated, so it is open in horizontal topology. 
 
 Set $ g(U_k(\alpha) \times \set{\beta}) \times \set{\root, \eps} $  is open since $ g $ is open.
 
 Consider point $(\pair{\alpha, k}, \beta)$, where $ k>0$;  it is isolated in the horizontal topology. It is easy to check that its $f$-image is also an isolated point.

 ($T^h_1$-continuity) If $ n>0 $, then points $ (\pair{\alpha, n}, \beta) $ and $ f(\pair{\alpha, n}, \beta) $ are both isolated in corresponding horizontal topologies. 
 
 Let us consider the case when $ n=0 $. 
 The minimal neighborhood of point $ f(\pair{\alpha, 0}, \beta) = \pair{g(\alpha, \beta), \root}$ equals
 \[ 
 R_1'(\pair{g(\alpha, \beta), \root}) =  R_1(g(\alpha, \beta)) \times \set{\root, \eps}.
 \] 
 We need to check that preimage of this neighborhood is open in horizontal topology of space $\X \times \Y$. To do this, we ensure that any point of this preimage has a neighborhood that is also in the preimage.

 Let $f(\pair{\gamma, k}, \delta) = \pair{g(\alpha, \beta), \eps}$. In this case $k>0$ and point $\pair{\pair{\gamma, k}, \delta}$ is isolated in $T_1^h$.
 
 Now let
 $$
 f(\pair{\gamma, k}, \delta) \in R_1(g(\alpha, \beta)) \times \set{\root}.
 $$ 
 Then $k=0$ and $g(\gamma, \delta) \in R_1(g(\alpha, \beta))$. Since $g$ is continuous, there exists a neighborhood $U_l(\gamma) \times \set{\delta}$, such that its $g$-image is in $R_1(g(\alpha, \beta))$. Hence, $f$-image of neighborhood $U'_m(\gamma, 0) \times  \set{\delta}$, where $m = \max(k,l) = l$ is included in $R_1(g(\gamma, \delta)) \times \set{\root, \eps} \subseteq R_1(g(\alpha, \beta)) \times \set{\root, \eps}$.
 
 	($T^v_2$-openness) Consider $T^v_2$-neighborhood of point $(\pair{\alpha, 0}, \beta)$  
 $$ \set{\pair{\alpha, 0}} \times U_k(\beta) .$$ 
 Its image is 
 \begin{align*}
 	f(\set{\pair{\alpha, 0}} \times U_k(\beta)) &=
 	\setdef[\pair{g(\alpha, \beta'), \root}]{\beta' \in U_k(\beta)} =\\
 	&= g(\set{\alpha}\times U_k(\beta)) \times \set{\root}.
 \end{align*}
 And this is a vertically opened set since $ g $ is a p-morphism.
 
Let $ n>0 $ then $ \set{\pair{\alpha, n}} \times U_k(\beta) $ is a vertical neighborhood of point $(\pair{\alpha, n}, \beta)$. Then
\begin{align*}
	f(\set{\pair{\alpha, n}} \times U_k(\beta)) &=
	\setdef[\pair{g(\alpha\lceil_n 0^\omega, \beta'\lceil_n 0^\omega), f_{\omega}(\gamma)}]{\beta' \in U_k(\beta), \beta' = \beta'\lceil_n\;\cdot\;\gamma} =\\
	&=\bigcup_{\beta'\in U_k(\beta)} \setdef[\pair{g(\alpha\lceil_n 0^\omega, \beta'\lceil_n 0^\omega), f_{\omega}(\gamma)}]{\beta' = \beta'\lceil_n\;\cdot\;\gamma}.
\end{align*} 
 
Note that $ U_k(\beta) $ can be represented as a union of similar neighborhoods but with greater $k$. Hence, we can prove vertical openness of $f$ for all $k$ greater than some number. Let $ k>\max (n, st(\alpha)) $ then $\alpha\lceil_k 0^\omega= \alpha$ and 
	\begin{align*}
	\bigcup_{\beta'\in U_k(\beta)}
	\setdef[\pair{g(\alpha\lceil_n 0^\omega, \beta'\lceil_n 0^\omega), f_{\omega}(\gamma)}]{\beta' = \beta'\lceil_n\;\cdot\;\gamma} &=\\
	= \setdef[\pair{g(\alpha, \beta\lceil_n 0^\omega), f_{\omega}(\gamma)}]{\beta'\in U_k(\beta)\ \beta' = \beta' \lceil_n\;\cdot\;\gamma} &= R'_2 (f(\pair{\alpha, k}, \beta)).
\end{align*}
Indeed, let $f(\pair{\alpha, n}, \beta) = \pair{\vec a, \vec b}$, where $\vec a = g(\alpha, \beta\lceil_n 0^\omega)$ and $\vec b = f_\omega (\gamma)$, $\beta = \beta \lceil_n \gamma$, then
\[ 
R'_2 \left (\pair{\vec a, \vec b}\right ) = \setdef[\pair{\vec a, \vec b \cdot \vec c}]{\vec c \in T_2}.
\]
Hence, $\beta' = \beta\lceil_n \cdot \vec c \cdot 0^\omega$ and $f(\pair{\alpha, n}, \beta') = \pair{\vec a, \vec b \cdot \vec c} $. The inverse is also true since for any $\gamma$
$ 
f(\pair{\alpha, n}, \beta\lceil_k \cdot \gamma) = \pair{\vec a, \vec b \cdot f_\omega(\gamma)} \in R'_2 \left (\pair{\vec a, \vec b}\right ).
$

($T^v_2$-continuity) Let us show that $f^{-1}\left ( R'_2 \left (f(\pair{\alpha, n}, \beta) \right ) \right )$ is vertically open.  
We assume that $f(\pair{\alpha, n}, \beta) = \pair{\vec a, \vec b}$ and $f\left (\pair{\alpha', m}, \beta'\right ) \in R'_2 \left (\pair{\vec a, \vec b}\right )$. 
Hence, 
\begin{align*}
	f\left (\pair{\alpha', m}, \beta'\right ) &= \pair{\vec a, \vec b \cdot \vec c},\\    
	g(\alpha'\lceil_m 0^\omega, \beta'\lceil_m 0^\omega) &= g(\alpha\lceil_n 0^\omega, \beta\lceil_n 0^\omega) = \vec{a} \in T_{2,2},\\
	f_\omega (\gamma') = \vec{b} \cdot \vec{c}&\hbox{ and }f_\omega (\gamma) = \vec{b}\hbox{, where }\beta' = \beta'\lceil_m \cdot \gamma'\hbox{ and }\beta = \beta\lceil_n \cdot \gamma.
\end{align*}

Let $k>\max(m,n, st(\alpha'), st(\beta'))$, then the $f$-image of neighborhood $\set{\alpha'} \times U_k(\beta')$ is included in $R'_2 \left (\pair{\vec a, \vec b}\right )$.
\end{proof}

Theorem \ref{thm:main} follows from Lemma \ref{lem:S41xS4corr}, Lemma \ref{lem:S41xS4comp} and Theorem \ref{thm:pmorphism}.

\section{Conclusions}
We are only in the beginning of the path in researching topological and neighborhood products. The neighborhood product is a generalization of the topological product for arbitrary modal logics. By now only very basic problems are solved. Further research can include the following problems:
\begin{enumerate}
	\item Find logics $ \logic{S4.1} \times_t \logic{S4.1} $, $ \logic{S4.2} \times_t \logic{S4} $, $ \logic{S4.2} \times_t \logic{S4.1} $ and $ \logic{S4.2} \times_t \logic{S4.2} $. 
	
	\item Find sufficient conditions for the topological product of logics to coincide with the fusion (the semiproduct or the product) of corresponding logics. 
	
	\item By interpreting modal operator $ \romb $ with the derivational operator in a topological space one can prove the topological completeness for logics weaker then   $ \logic{S4} $ (see \cite{bezhanishvili2005some}).
	In \cite{kudinov2018IGPL,kudinov_aiml12} it was proved that $\logic{D4} \times_t \logic{D4} = \logic{D4} \ast \logic{D4}$ and $\logic{K4} \times_t \logic{K4} = \logic{K4} \ast \logic{K4} + \Delta$, where $ \Delta $ is a certain set of closed formulas.
	In this setting one can study the products of extensions of logics $ \logic{D4} $ and $ \logic{K4} $ with axioms $A1$ and $A2$.
\end{enumerate}


\bibliographystyle{aiml22}
\bibliography{kudinov}

\end{document}